\documentclass[11pt,final]{amsart}
\usepackage{amsmath,amssymb,amsthm,showkeys}

\usepackage{amsmath}
\usepackage{graphicx}
\usepackage{color}

\newtheorem{thm}{Theorem}[section]
\newtheorem{lem}{Lemma}[section]
\newtheorem{corollary}{Corollary}[section]
\newtheorem{prop}{Proposition}[section]

\newtheorem{claim}{Claim}

\theoremstyle{definition}

\theoremstyle{remark}
\newtheorem{remark}{Remark}

\numberwithin{equation}{section}

\newenvironment{Abstract}
{\begin{center}\textbf{\footnotesize{Abstract}}%
\end{center} \begin{quote}\begin{footnotesize}}
{\end{footnotesize}\end{quote}\bigskip}

\newcommand{\ga}{\gamma}

\newcommand{\SSS}{\mathcal{S}}

\newcommand{\abs}[1]{\left | #1\right |}
\newcommand{\round}[1]{{\left ( #1 \right )}}
\newcommand{\norm}[2]{{\left\| #1 \right\|}_{#2}}
\newcommand{\f}[2]{\frac{#1}{#2}}

\newcommand{\set}[1]{\{ #1\}}
\newcommand{\jnorm}[1]{\langle#1\rangle}

\newcommand{\oop}{\frac{1}{p}}

\newcommand{\ooq}{\frac{1}{q}}

\newcommand{\ve}{\varepsilon}

\newcommand{\vp}{\varphi}

\newcommand{\rn}{{\mathbf R}^n}

\newcommand{\rone}{\mathbf R}

\newcommand{\nn}{\mathbf N}

\newcommand{\zz}{\mathbf Z}

\newcommand{\lp}{L^{p}}

\newcommand{\lqq}{L^{q}}
\newcommand{\lr}{L^{r}}

\newcommand{\lo}{L^{1}}
\newcommand{\lt}{L^{2}}

\newcommand{\lmixpq}{L^p_tL^q_x}

\newcommand{\ch}{\mathcal H}

\newcommand{\intl}{\int\limits}

\newcommand{\suml}{\sum\limits}
\newcommand{\supl}{\sup\limits}
\newcommand{\supp}{{\rm supp \, }}

\newcommand{\p}{\partial}

\newcommand{\beq}{\begin{equation}}
\newcommand{\eeq}{\end{equation}}
\newcommand{\beqna}{\begin{eqnarray*}}
\newcommand{\eeqna}{\end{eqnarray*}}
\newcommand{\beqn}{\begin{equation*}}
\newcommand{\eeqn}{\end{equation*}}
\newcommand{\bep}{\begin{proof}}
\newcommand{\enp}{\end{proof}}
\newcommand{\bprop}{\begin{proposition}}
\newcommand{\eprop}{\end{proposition}}
\newcommand{\bt}{\begin{theorem}}
\newcommand{\et}{\end{theorem}}
\newcommand{\bex}{\begin{Example}}
\newcommand{\eex}{\end{Example}}
\newcommand{\bc}{\begin{corollary}}
\newcommand{\ec}{\end{corollary}}
\newcommand{\bcl}{\begin{claim}}
\newcommand{\ecl}{\end{claim}}
\newcommand{\bl}{\begin{lemma}}
\newcommand{\el}{\end{lemma}}

\thanks{M.T.
 is partially
supported by Contract PRIN-2005010482 (2005-2007) with Italian
Ministry of Education, University and Research and the INDAM
Contract ``Mathematical modeling and numerical analysis of quantum
systems with applications to nanosciences''.}

\begin{document}

\title[Commutator Estimate]{On a Calder\'on-Zygmund commutator-type
estimate}

\author{Mirko Tarulli}

\address{Mirko Tarulli\\
Dipartimento di Matematica Universit\`a di Pisa\\
Largo B. Pontecorvo 5, 56100 Pisa, Italy}

\email{tarulli@mail.dm.unipi.it}

\author{J. Michael Wilson}

\address{J. Michael Wilson\\
Mathematics Department University of Vermont 
16 Colchester Ave. Burlington, Vermont 05405, Usa }

\email{wilson@cems.uvm.edu}





\maketitle

\date{}

\begin{Abstract}
In this paper we extend a Calder\'on-Zygmund commutator-type
estimate. This estimate enables us to prove an embedding result
concerning weighted function spaces.
\end{Abstract}

{\bf 2000 Mathematics Subject Classification: 42B15, 42B20, 42B25,
\\42B35.}
\vskip+0.2cm

{\bf Keywords and phrases:} Calder\'on-Zygmund commutators, Hardy-Littlewood Maximal Operator, Littlewood-Paley Theory, Weighted Besov Spaces. \vskip+0.5cm

\section{Introduction and statement of results}

The theory of commutators has been extensively studied. They were
introduced by R. Coifman, R. Rochberg and G. Weiss in the paper
\cite{CRW}, in which it is proved that, given a singular integral
$T$ with standard kernel (we refer to \cite{Du} and \cite{Stein1}
for more details), the operator $[b,T]=bT-Tb$ is bounded in $L^p,\,
1<p<\infty$ if b is a $BMO$ function; the converse implication is
due to S. Janson, who stated it in \cite{Jans}. We refer also to the
remarkable paper of C. P\'erez \cite{Perez}, in which he underlined
that $[b,T]$ is more singular than the operator $T$, proving that
this commutator does not satisfy the corresponding $(1,1)$ estimate,
but a weaker one, and to the paper of A. Uchiyama \cite{Uch} in
which is proved that the commutator of $T$ and $b$ is a compact
operator provided that $b\in VMO$.  We also have a more general form
of commutator operators; we
refer to the papers of A. P. Calder\'on \cite{Cal1} and \cite{Cal2},
in which the author proved that these kinds of operators are bounded
on $\lt$ and, under certain condition, they are of
Calder\'on-Zygmund-type (see also his survey paper \cite{Cal3} and
to the paper of R. Coifman and Y. Meyer \cite{CM}).

In this paper we generalize a Calder\'on-Zygmund commutator-type
estimate introduced by T. Tao in the paper \cite{tao1} and
successively presented in a more extended form by T. Tao and I.
Rodnianski in the paper \cite{rodntao} which reads
\begin{equation}\label{eq.RT1}
\norm{[P_k,f] g}{L^r}\leq C 2^{-k} \norm{\nabla
f}{L^q}\norm{g}{L^p},
\end{equation}
for any Schwartz functions $f,g$ whenever $1\leq r,p,q\leq \infty,$
$1/r=1/q+1/p$ and with $[P_k,f]g=P_k(fg)- fP_k(g).$ Here $P_k $ is
the $k^{th}$ Littlewood-Paley projection (see Section \ref{NotPre}
for the definition).

We prove that this estimate can be improved by considering not only
derivatives of first order---the operator $\nabla=(\p_{x_1}, \ldots
\p_{x_n})$ appearing on the right of \eqref{eq.RT1}---but also
introducing a homogeneous differential operator of fractional order
$s$ belonging to a bounded interval of $\rone.$ This is connected to
the fact that the structure of the commutator operator is very rich,
and that we can use the good properties of the kernels of such
operators coupled with some cancelation involving terms in its
integral representation. On the other hand, we can also exploit the
relation between differential operators of negative fractional order
and $L^p$ ($1<p<\infty$) bounded singular integral operators with
homogeneous kernels.

Our main result in this direction is:
\begin{thm}\label{CommThm}
Given a function $P\in\lo(\rn)$ which satisfies, for $\ve>0$,
$0<s<1$, and $j\in \nn \setminus \{0\} ,$
\begin{eqnarray}\label{eq.MainAss}
\intl_{2^{j-1}\leq \abs{x}\leq 2^{j}}\abs{P(x)} \, dx \lesssim 2^{-j(\ve+s)}, \ \ \ \ \   \intl_{\abs{x}\leq 1}\abs{P(x)}\, dx\lesssim 1,
\end{eqnarray}
if we set $P_k=2^{nk}P(2^{k}x)$, for all $k\in \zz$, we have, for
any $f,\,g\in \SSS(\rn)$,
\begin{eqnarray}\label{eq.Main}
\norm{\int P_k(y)\round{f(x)-f(y)}g(x-y) \, dy }{L^{r}}\lesssim 2^{-ks}\norm{\abs{D}^s f}{\lp}\norm{g}{\lqq},
\end{eqnarray}
provided that $1\leq p,q,r \leq \infty, \ \oop+\ooq=\f1r$.
\end{thm}

The theorem implies:

\bc \label{CaldZygComm}
The estimate
\begin{equation}\label{eq.TW}
\norm{[P_k,f] g}{L^r}\leq C 2^{-sk} \norm{\abs{D}^s f}{L^p}\norm{g}{L^q},
\end{equation}
holds for all $k\in \zz,$ all functions $f, g\in \SSS(\rn)$ and all
$1\leq p,q,r \leq \infty$ such that $\oop+\ooq=\f1r$ and $0<s<1.$
\ec

We apply these commutator estimates to problems of well-posedness
and local and global existence for a class of evolution equations.
In particular, we are interested in certain function spaces that
come up in obtaining the relevant \emph{a priori} estimates for
these equations. These equations can be written in the general form,

      \begin{equation}\label{DisrEqF}
        \begin{cases}
    i u_t(t,x) + P(D) u(t,x) = F \ , \ \ t \geq 0 \ , \ x \in \rn\\
         \\
  u(0,x) = u_0(x),

   \end{cases}
    \end{equation}
where the symbol $P(\xi)$ (it could be either a scalar or a real
matrix) has some regularity properties. Some of the most important
\emph{a priori} estimates on solutions of problems of this type have
the form
\begin{equation}\label{apriori1}
\norm{\abs{D_x}^\alpha v}{L^2_tB_{w}'}\lesssim \norm{f}{\ch}+\norm{\abs{D_x}^\beta F}{L^2_tB_{w}},
\end{equation}
with $\alpha, \beta \in \rone$ satisfying some relations, $B_w$ is a
suitable weighted Banach space (and $B_w'$ is its dual space), and
$\ch$ is a Hilbert space. In fact, for the case of the Schr\"odinger
equation (i.e., with $P(D)=-\Delta$ in \eqref{DisrEqF}), C. Kenig,
G. Ponce and L. Vega established the so-called smoothing estimates
in \cite{KPV93} (see also the papers of A. Ruiz and L. Vega
\cite{RuVe94}, B. Perthame and L. Vega \cite{PV}, and V. Georgiev
and M. Tarulli \cite{GeTa}). These were later extended to more
general second order Schr\"odinger equations in \cite{KPV98}. The
smoothing estimate reads:
\begin{eqnarray}\label{eq.3214}
   \supl_m 2^{-m/2} \norm{|D_x|^{1/2} u}{L^2_t L^2(|x|\sim
2^{m})}&\leq& C\norm{u_0}{L^2}+\nonumber\\
&+&\suml_m 2^{m}  \norm{|D_x|^{-1/2} F}{L^2_t L^2(|x|\sim
2^{m})},
\end{eqnarray}
and corresponds to \eqref{apriori1} if one chooses $\alpha=1/2,$
$\beta=-1/2$ and $B_w=L^1(\zz, 2^{ks}dk, L^2(|x|\sim2^m)).$ Chihara
proved in \cite{Chi} that similar estimates can be obtained for the
solutions of a larger class of partial differential equations of the
form \eqref{apriori1}. More precisely, if we take $\alpha=d/2,$
$\beta=-d/2$ and $B_w=L^2(\rn, \jnorm{x}^{\delta/2}dx):=L^2_\delta $
the estimates are
\begin{eqnarray}\label{eq.Chi}
    \norm{|D_x|^{d/2} u}{L^2_t L^2_{-\delta}}\lesssim \norm{u_0}{L^2}+ \norm{|D_x|^{-d/2} F}{L^2_t L^2_{\delta}},
\end{eqnarray}
with $\jnorm{x}=\sqrt{1+\abs{x}^2},\, d>1, \delta >0,$ where the
assumptions on $P(D)$ reflect the dispersive nature of the equation
\eqref{DisrEqF}; or in other words, its symbol $P(\xi)$ behaves like
$\abs{\xi}^m,\, m>1.$ We point out that this estimate is also valid
if one uses some rotation invariant norms, just seen for the case of
the Schr\"odinger equation, and which give scale-invariant
estimates. There are analogous estimates, transposed to the
framework of the wave equation, which correspond to the case $ P(D)=
\begin{pmatrix}
 0 & -I \\
 -\Delta & 0
 \end{pmatrix}
, \, u_0=(f,g),$ in \eqref{DisrEqF}, as showed by N. Burq in
\cite{Brq} and M. Tarulli in \cite{Ta}. These estimates have the
form,
\begin{eqnarray}\label{eq.Wave}
    \norm{|D_x|^{\ga} \round{u, |D_x|^{-1} \p_t u}}{L^2_t(\rone, (L^2_{-\delta}\times L^2_{-\delta})}&\lesssim& \norm{(f,g)}{\dot{H}^{\ga}\times \dot{H}^{\ga-1} }+\nonumber\\
&+&  \norm{|D_x|^{\ga-1} F}{L^2_t L^2_{\delta}},
\end{eqnarray}
which are valid as well if one uses norms that give scale-invariant
estimates, instead of weighted Sobolev norms.

On the other hand, many authors (see, e.g., \cite{rodntao} or
\cite{Ho83}) underline the necessity of using Besov (weighted)
spaces instead of Sobolev (weighted) ones, because they allow one to
better control some nonlinear terms arising in the evolution problem
\eqref{DisrEqF}. (An example is the case of the magnetic potential
perturbation $F=-A(x)\cdot \nabla u,$ with some decay property
assumed on the function $A(x)$; for more detail we refer again to
\cite{GeSteTa}). The Besov norms give not only natural
scale-invariant estimates, important for many problems in partial
differential equations, but also enable us to exploit the
Littlewood-Paley theory on which they are built. There is a vast
literature on this subject. We have found many results on
equivalences of Besov norms and Sobolev norms (we refer to
\cite{BL}), but very few regarding the relation between weighted
Besov norms and weighted Sobolev norms. A first attempt to do this
was due to V. Georgiev, A. Stefanov and M. Tarulli \cite{GeSteTa};
we generalize those results here. Motivated by this, we shall show
that weighted $L^2$-based Sobolev spaces, defined as the closure of
Schwartz functions $\phi$ with respect to the norms
\begin{equation}\label{eq.weiSob}
\sum_m 2^{m/2} \norm{D_x^{\ga-d/2} \phi(x)}{L^2(|x|\sim
2^m)}
\end{equation}
are embedded in suitable weighted Besov spaces $Y^{\ga,d}$,
defined as the closure of the functions
$$\phi(x) \in C_0^\infty( \rn \setminus
 \{ 0\})$$ with respect to the norm
\begin{equation}\label{eq.defY}
   \norm{\phi}{Y^{\ga,d}}:=\left(\suml_k 2^{2\ga k}  \norm{\phi}{Y_{d,k}}^2\right)^{1/2}.
\end{equation}
In \eqref{eq.defY} the spaces $Y_{d,k}$ are defined by the
 norms\footnote{As noted in \cite{GeSteTa}, the expressions $\phi\to \norm{\phi}{Y_{d,k}}$
are not faithfull norms, in the sense that may be zero, even for
some  $\phi\neq 0$. On the other hand, they satisfy all the other
norm requirements and $\phi\to (\sum_k 2^{\ga k}\|\phi_k\|_{Y^{d,k}}^2)^{1/2}$
is a norm!}
\begin{equation}\label{eq.defY2}
\norm{\phi}{Y_{d,k}}= 2^{- dk/2} \sum_m 2^{m/2} \norm{\vp(2^{-m}\cdot)\phi_k}{
L^2_x},
\end{equation}
where $\phi_k=P_k\phi$ (see again section \ref{NotPre}).
Similarly, we have the Banach space $\overline{Y}^{d,\ga}$, with
norm
\begin{equation}\label{eq.defYdual}
   \norm{\phi}{\overline{Y}^{d,\ga}}:=\left(\suml_k 2^{2\ga k}  \norm{\phi}{Y'_{d,k}}^2\right)^{1/2},
\end{equation}
where  the spaces $Y'_{d,k},$ equipped with the norm
\begin{equation}\label{eq.defYdual2}
\norm{\phi}{Y'_{d,k}}= 2^{dk/2} \sup_m 2^{-m/2} \norm{\vp(2^{-m}\cdot)\phi_k}{
L^2_x},
\end{equation}
are duals to the spaces $Y_{d,k}.$

From the fact that  $\|F\|_{L^2(|x|\sim 2^m)} \sim
\|\vp(2^{-m}\cdot) F\|_{L^2}$ we can replace at will,
 $\|F\|_{L^2(|x|\sim 2^m)}$, with the comparable expression,  $ \|\vp(2^{-m}\cdot) F\|_{L^2}$.
 We will often do this in the sequel, in order to make use of localizations in both the space and
 frequency variables.

We shall show:

\begin{thm}
\label{th:4prima} Assume that $n\geq2.$ There is a constant $C=C(n)$, so that for every
function $\phi \in \SSS( \rn)$ and $\ga, d\geq 0$, we
have
\begin{equation}
\label{eq:09prima}\norm{\phi}{Y^{\ga,d}} \leq C  \sum_m 2^{m/2}
\norm{ \vp(2^{-m}\cdot)D_x^{\ga-d/2} \phi(t,x)}{L_x^2}.
\end{equation}
\end{thm}

We conclude this section proving the following corollary.

 \bc
\label{le:4} Assume that $n\geq2.$  There is a constant $C=C(n)$, so that for every
function $\phi \in \SSS( \rn)$ and $\ga, d\geq 0$, we
have
\begin{equation}
\label{eq:09mod}
\supl_m 2^{-m/2} \norm{\vp(2^{-m}\cdot)D_x^{\ga+d/2} \phi(t,x)}{
L^2_x} \leq C_n \norm{\phi}{\overline{Y}^{\ga,d}}.
\end{equation}
\ec

\bep This corollary is an easy consequence of the previous one. We
can give its proof in one line. Consider the mapping $\phi\to
D_x^{-\ga} \phi$, and take the estimate dual to \eqref{eq:09prima}.
That proves \eqref{eq:09mod}. \enp

\begin{remark}
We see that we can also assume that the functions $\phi$
are dependent on the time variable, that is $\phi=\phi(t,x)$. We just 
replace in \eqref{eq.weiSob}, \eqref{eq.defY2} and
\eqref{eq.defYdual2} the $L^2_x$ norm by the mixed one $L_t^2L^2_x,$ this doesn't
cause confusion, because the whole proof only depends on what happens in the space variable
(see, fore more details, the proofs of Lemma \ref{WBS} and Theorem \eqref{th:4prima}. 
The estimates \eqref{eq:09prima} and \eqref{eq:09mod} remain valid and
defined in this way, the spaces $Y^{\ga,d}$ (now based on $L_t^2L^2_x$)  are related to the others
listed above and connected with the evolution problem
\eqref{DisrEqF}. If we choose $\ga=0,$ we are in the case of
\eqref{eq.Chi}. Furthermore, the smoothing spaces in \eqref{eq.3214}
are a particular case of the previous one; that is. $\ga=0$ and
$d=1$ while the others in \eqref{eq.Chi} are derived if we pick $\ga=0$ and
$d=\tilde{d}$. Finally, for the wave equation, we may consider
$\ga=\tilde{\ga}-1/2$ and $d=1.$ The authors say that it is possible
to use these results in many other situations.
\end{remark}

\begin{remark}
Theorem \eqref{th:4prima} and Corollary \eqref{le:4} generalize completely the estimates proved in Lemma 1 in
the paper \cite{GeSteTa}, which correspond to \eqref{eq:09prima} and \eqref{eq:09mod}
in the case $\ga=0$ and
$d=1.$
\end{remark}
\subsection{Outline of paper.} Section \ref{NotPre} is devoted to
presenting the ancillary theory used in the whole paper. In Section
\ref{CommMain} we prove Theorem \ref{CommThm}, using two important
introductory lemmas. In Section \ref{Emb}, we give the
proof of the embedding Theorem \ref{th:4prima}. Finally, in the Appendix
we present an alternative proof of Lemma \ref{CommBump}.

\vspace{0.6cm} {\it Acknowledgement:} The authors are grateful to
Atanas Stefanov for an interesting discussion regarding especially
the proof of the second part of the Theorem \ref{th:4prima}.

\section{Notations and preliminaires}\label{NotPre}

In this section we introduce some important tools and useful notation

\vspace{0.3cm}
$\mathbf{1.}$ Given any two positive real number $a, b,$ we write $a\lesssim b$ to indicate
$a\leq C b,$ with $C>0.$ In the same way we say that $a\sim b,$ if there
exist poitive constants $c_1,  c_2,$ so that $c_1 a\leq b \leq c_2 a.$

\vspace{0.3cm}
$\mathbf{2.}$
The Fourier transform  and its inverse are defined by
\begin{eqnarray*}
& &
\hat{f}(\xi)= \intl_{\rn} f(x)e^{- 2\pi i x\cdot \xi } dx\\
& &
 f(x)=\intl_{\rn} \hat{f}(\xi)e^{2\pi i  x\cdot \xi } d\xi.
\end{eqnarray*}
$\textbf{3.}$ Consider a positive, decreasing, function
$\psi:\rone_+ \to \rone_+ $, which is smooth away from zero,
supported in
 $\{\xi:0\leq \xi \leq 2\}$ and equals $1$ for all $ \ 0 \leq \xi\leq
 1$. Set $\vp(\xi)\equiv\psi(\xi)-\psi(2\xi)$. The function $\vp$ is
 smooth and satisfies:
\begin{eqnarray}\label{pallit}
 \supp \vp(\xi)  \subseteq \set{1/2\leq |\xi|\leq 2}, \ \ \,
  \sum_{k\in \zz} \vp(2^{-k} \xi)=1, \, \, \, \,  \forall \ \xi \neq 0
\end{eqnarray}
We use $\vp(\xi)$ to create a \emph{frequency space localization}. Similarly, we can introduce the notion of \emph{physical space
localization} by means of the function $\vp(x)$ which satisfies
\eqref{pallit}. In higher dimensions, we slightly abuse notation and
denote a function with similar properties by the same symbol;
namely, $\vp(\xi)=\vp(|\xi|)$, $\psi(x)=\psi(|x|)$ etc. Note that
for $n>1$, $\psi:\rn\to\rone$ is a smooth function even at zero. The
$k^{th}$ Littlewood-Paley projection operator $P_k$, $P _{\leq k}$
is defined by
$$\widehat{P_k f}(\xi)=\vp(2^{-k}\xi) \hat{f}(\xi),$$ and $$\widehat{P _{\leq k} f}(\xi)=\psi(2^{-k}\xi) \hat{f}(\xi),$$ we also write $P_{<k}$ for $P_{\leq k-1}.$
Let us observe that the kernel of $P_k$ is integrable, smooth and real-valued
for every $k$ and that every $P_k$ is bounded on every
$L^p$ space, $1\leq p\leq \infty$, and that it commutes with
constant coefficient differential operators. We may also consider
$P_{>k}:=\sum_{l>k} P_l$, which essentially restricts the Fourier
transform to frequencies $\gtrsim2^k$. In the paper we shall use,
without any discussion, the notation $f_k, \ f_{\leq k}, \ f_{<k}$,
to indicate respectively $P_k f, \ P _{\leq k} f, \ P _{<k} f,$ and,
for any fixed $N\geq 1,$ we usually write $f_{k-N \leq \cdot \leq
K+N}$ for $P _{\leq k+N} f- P _{<k-N} f.$ An obvious modification give also the operators
 $P _{< k+N} f- P _{\leq k-N} f.$ and consequently we have $f_{k-N < \cdot <
K+N}.$ We notice also that one
can obtain the representations
\begin{equation}\label{eq.PLint1}
P_k f(x)=2^{nk}\int \widehat{\vp}(2^{k}y) f(x-y)dy
\end{equation}
and
\begin{equation}\label{eq.PLint2}
P_{\leq k} f(x)=2^{nk}\int \widehat{\psi}(2^{k}y) f(x-y)dy,
\end{equation}
for any $f\in \SSS(\rn).$ 
Another helpful observation is that, for
the differential operator $D_x^s$, defined via
 the multiplier $|\xi|^s$, we have the identity
\beq\label{eq.tilLP}
D_x^s P_k u=2^{ks} \tilde{P}_k u,
\eeq
 where
$\tilde{\vp}(\xi)= \vp(\xi)|\xi|^s$. Notice that the operator $\tilde{P}_k$ has kernel satisfying the same estimate of $P_k.$ In view of this we say that the two operators are equivalent
and we make no difference in the use of them.

\vspace{0.3cm} $\mathbf{5.}$ The Littlewwod-Paley theory just presented enables us to
introduce the Bernstein inequality 
for $\rn.$ Given a function $f\in \SSS(\rn),$ for any dyadic number $K$
the estimate
\beq\label{eq.Bern}
\norm{P_{\leq K}f}{\lqq}\lesssim K^{n\round{\frac{1}{p}-\frac{1}{q}}}\norm{P_{\leq K} f}{\lp}
\eeq
is valid provided that $1\leq p\leq q\leq \infty.$ This inequality remains true
if we replace the operator $P_{\leq K}$ by the other one $P_{K}$ introduced above.

\vspace{0.3cm} $\mathbf{6.}$ We say a function $f$ belongs to
$\lr(X)$ if the norm
$$\| f \|_{\lr_x(X)} =\left(
\int_{X} |f(x)|^r dx\right)^{1/r}$$ is finite; we will usually omit
the domain of integration if $X=\rn$. In a similar way we say a
function $u$ belongs to $\lmixpq$ if the mixed space-time norm,
$$
\| u\|_{L^q_t L^{r}_x} = \left(\int_{\rone} \left(
\int_{\rn} |u(t,x)|^r dx\right)^{q/r} dt \right)^{1/q},
$$is finite.

\vspace{0.3cm} $\mathbf{7.}$
The Hardy-Littlewood maximal function is defined by
\begin{equation}\label{eq.Hardy-Littl}
M(f)(x)=\supl_{r>0}\frac{1}{\mu\round{{B(x,r)}}}\int_{B(x,r)}\abs{f(y)}dy,
\end{equation}
where $B(x,r)$ is the ball of radius $r$ centered at $x$ and $\mu$
is the classical Lebesgue measure on $\rn.$ We remember also that
the operator $M$ is bounded on $L^{p}$ for $1<p\leq \infty$ (we
refer to \cite{Du} or \cite{Stein1} for more details).

\section{Proof of the Theorem \ref{CommThm}}\label{CommMain}

First of all we need the following lemma

\begin{lem}\label{localHo}
Given any function $f\in \SSS(\rn)$ and real number $0<s<1,$ there
exists a positive constant $C_{s,n}$ so that, for all $x$ and $y$ in
$\rn$,
\begin{equation}
\label{eq.locHol}
\abs{f(x)-f(y)}\leq C_{s,n} \abs{x-y}^s \round{M(\abs{D}^s f)(x)+M(\abs{D}^s f)(y)},
\end{equation}
where $M(\cdot)$ is the Hardy-Littlewood maximal operator.
\end{lem}

\begin{proof}
If $f\in \SSS(\rn)$, we can easily calculate its inverse Fourier
transform,  obtaining
\begin{eqnarray}\label{eq.RieszPot}
f(x)=\int \hat{f}(\xi)e^{2\pi i \langle x,\xi \rangle} d\xi=\intl
\abs{\xi}^{-s}  \abs{\xi}^{s} \hat{f}(\xi)e^{2\pi i \langle x,\xi
\rangle} d\xi.
\end{eqnarray}
If $0<\alpha<n$,
$$ \abs{\xi}^{-\alpha}=c_{\alpha,n}\widehat{\abs{x}^{\alpha-n}}(\xi).$$
Therefore, we can rewrite the right-hand side of \eqref{eq.RieszPot}
as
\begin{eqnarray}\label{eq.RieszPot2}
\int \abs{\xi}^{-s} \round{ \abs{\xi}^{s} \hat{f}(\xi)}e^{2\pi i \langle x,\xi \rangle} d\xi=c_{s,n}\int \frac{\abs{D}^s f(t)}{\abs{x-t}^{n-s}} dt.
\end{eqnarray}
Combining \eqref{eq.RieszPot} and  \eqref{eq.RieszPot3},  we have
the representation\footnote{The kernel in the integral
\eqref{eq.RieszPot3} is in generally referred to as a Riesz
potential.};
\begin{eqnarray}\label{eq.RieszPot3}
f(x)=c_{s,n}\int \frac{\abs{D}^s f(t)}{\abs{x-t}^{n-s}} dy,
\end{eqnarray}
so, we can write:
\begin{eqnarray}\label{eq.RieszPot4}
f(x)-f(y)&=&c_{s,n}\int \frac{\abs{D}^s f(t)}{\abs{x-t}^{n-s}} dt-c_{s,n}\int \frac{\abs{D}^s f(t)}{\abs{y-t}^{n-s}} dt=
\nonumber\\
&=&c_{s,n}\int \abs{D}^s f(y) \round{\frac{1}{\abs{x-t}^{n-s}} - \frac{1}{\abs{t-y}^{n-s}}} dt
\nonumber\\
&=&I+II,
\end{eqnarray}
where
\begin{equation}\label{eq.RieszPot5a}
I=c_{s,n}\intl_{\abs{x-t}\leq 3 \abs{x-y}} \abs{D}^s f(t) \round{\frac{1}{\abs{x-t}^{n-s}} - \frac{1}{\abs{y-t}^{n-s}}} dt,
\end{equation}
and
\begin{equation}\label{eq.RieszPot5b}
II=c_{s,n}\intl_{\abs{x-t}\geq 3 \abs{x-y}} \abs{D}^s f(t) \round{\frac{1}{\abs{x-t}^{n-s}} - \frac{1}{\abs{y-t}^{n-s}}} dt.
\end{equation}
We treat $I$ first. We notice that it is easy to obtain the
following set inclusions:
$$\set{y: \abs{x-t}\leq 3 \abs{x-y}}\subseteq \set{y: \abs{t-y}\leq 4 \abs{x-y}}.$$
We can write:

\begin{eqnarray}\label{eq.RieszPot6}
 \abs{I}\leq c_{s,n}\round{\intl_{\abs{x-t}\leq 3 \abs{x-y}}\frac{ \abs{\abs{D}^s f(t) }}{\abs{x-t}^{n-s}} + \intl_{\abs{t-y}\leq 4 \abs{x-y}}\frac{\abs{\abs{D}^s f(t) }}{\abs{y-t}^{n-s}} dt}.
\nonumber
\end{eqnarray}\\
Given any $R\in \rone_+$ and $0<s<1$, we introduce the following
function:
\begin{equation}\label{def.Funct1}
\Phi_{1}^{R}(z)=\frac{1}{\abs{z}^{n-s}}\chi_{\set{\abs{z}\leq 4 R}};
\end{equation}
where $\chi_{\set{\abs{z}\leq 4 R}}$ is the characteristic function
of the set $\set{\abs{z}\leq 4 R}.$ We need to use the pointwise estimate
given in the following proposition (see section 2.2, Chapter III, \cite{Stein1}, or section 4, Chapter
II, \cite{Du}),
\begin{prop}\label{Stein}
if $\phi$ is a positive, radial decreasing function (on $(0, \infty))$ and integrable, then we have
$$\supl_{R>0}\abs{\Phi^{R}*g(z)}\leq
\norm{\Phi^{R}}\lo M(g)(z).$$
\end{prop}
In our case, from this lemma we get\footnote{Notice that this pointwise estimate is 
a weaker form of the one proved in \cite{Stein1} or in \cite{Du}.} $\abs{\Phi_{1}^{R}*g(z)}\leq
\norm{\Phi_{1}^{R}}\lo M(g)(z)$; i.e.,
\begin{equation}\label{eq.Stei1}
\intl_{\abs{z-t}\leq R} \frac{\abs{g(t)}}{\abs{z-t}^{n-s}} dt \leq
\norm{\Phi_{1}^{R}}\lo M(g)(z),
\end{equation}
for all $g\in \SSS(\rn).$ Then,
 coupling the estimates  \eqref{eq.RieszPot6},  \eqref{eq.Stei1} with $g=\abs{D}^s f $ and taking $\abs{x-y}\sim R,$ we obtain the inequalities

\begin{eqnarray}\label{eq.RieszPotFin1}
\abs{I}&\leq &C_{s,n} \round{\norm{\Phi_{1}^{R}}\lo M(\abs{D}^s
f)(x) +\norm{\Phi_{1}^{R}}\lo M(\abs{D}^s f)(y)}\leq
\nonumber\\
&\leq &C_{s,n}R^{s} \round{ M(\abs{D}^s f)(x) + M(\abs{D}^s f)(y) }.
\end{eqnarray}\\
where the third inequality in the above estimate
\eqref{eq.RieszPotFin1} is achieved using the fact that
$R^{s}=c_n\norm{\Phi_{1}^{R}}\lo.$

\vspace{0.2cm}
Now we look at $II.$ We have from \eqref{eq.RieszPot5b} the following estimate

\begin{equation}\label{eq.RieszPot7}
\abs{II}\leq c_{s,n}\intl_{\abs{x-t}\geq 3 \abs{x-y}}\abs{ \abs{D}^s f(t)} \abs{\frac{1}{\abs{x-t}^{n-s}} - \frac{1}{\abs{y-t}^{n-s}}} dt.
\end{equation}\\
We notice now that if $\abs{x-y}\geq 3 \abs{x-t},$ then the triangle
inequality easily yields
\begin{eqnarray*}
   2 \abs{x-y} \leq \abs{x-t}, \ \ \ \  \abs{y-t} \leq 2 \abs{x-t}, \ \ \ \  2 \abs{x-t} \leq 3\abs{y-t},
\end{eqnarray*}
which amount to saying that $\abs{x-t}\sim \abs{y-t}.$ This enables
us to obtain the following pointwise estimates:
\begin{eqnarray}\label{eq.RieszPot8}
 \abs{\frac{1}{\abs{x-t}^{n-s}} - \frac{1}{\abs{y-t}^{n-s}}}&\lesssim& \frac{\abs{x-t}^{n-s}-\abs{y-t}^{2n-2s}}{\abs{x-t}^{n-s}}\lesssim
 \nonumber\\
&\lesssim &\abs{x-y} \frac{\abs{x-t}^{n-s-1}}{\abs{x-t}^{2n-2s}}.
\end{eqnarray}
Replacing the esimate  \eqref{eq.RieszPot8} in the inequality \eqref{eq.RieszPot7}, we arrive at
\begin{equation}\label{eq.RieszPot9}
\abs{II}\leq c_{s,n} \abs{x-y} \intl_{\abs{x-t}\geq 3 \abs{x-y}}\frac{\abs{D}^s f(t)}{\abs{x-t}^{n-s+1}} dt.
\end{equation}
Let us introduce, for $R\in \rone_+$ and $0<s<1$, the following function
\begin{equation}\label{def.Funct2}
\Phi_{2}^{R}(z)=\frac{1}{\abs{z}^{n-s+1}}\chi_{\set{\abs{z}\geq 4 R}},
\end{equation}
where $\chi_{\set{\abs{z}\geq 4 R}}$ is the characteristic function of the set
$\set{\abs{z}\leq 4 R},$ by the same argument used to prove the inequality \eqref{eq.Stei1}, we see that the right hand side of the \eqref{eq.RieszPot9} satisfies the following
 \begin{equation}\label{eq.Stei2}
\intl_{\abs{z-t}\geq R} \frac{\abs{\abs{D}^s f(t)}}{\abs{z-t}^{n-s}} dt \leq \norm{\Phi_{2}^{R}}\lo M(\abs{D}^s f)(z).
\end{equation}
We thus have the following chain of inequalities

\begin{eqnarray}\label{eq.RieszPotFin2}
\abs{II}&\leq &C_{s,n} R\round{\norm{\Phi_{2}^{R}}\lo M(\abs{D}^s f)(x) +\norm{\Phi_{2}^{R}}\lo M(\abs{D}^s f)(y)}\leq
\nonumber\\
&\leq &C_{s,n}R^{s} \round{ M(\abs{D}^s f)(x) + M(\abs{D}^s f)(y) }.
\end{eqnarray}\\
 with $R^{s-1}=c_n\norm{\Phi_{2}^{R}}\lo$ and $\abs{x-y}\sim R.$
Finally, coupling the estimates \eqref{eq.RieszPotFin1} and
\eqref{eq.RieszPotFin2}, we obtain the bound \eqref{eq.locHol} as
desired.
\end{proof}

\begin{remark}
We notice that this lemma can be read as \lq\lq local H\"older
continuity.'' Up to this point, we have considered functions that
are small at infinity. It is quite natural to ask what happens for
more general ones. Given a function $f$ such that $\abs{D}^s f\in
L^\infty$, we obtain from \eqref{eq.locHol} the following:
\beq\label{eq.HolNorm}
 \supl_{x\neq
y}\frac{\abs{f(x)-f(y)}}{\abs{x-y}^s}<\infty. 
\eeq This means that
the function is H\"older continuous of order $0<s<1.$ Furthermore,
the H\"older continuous norm is related to another one, involving
the frequency decomposition \eqref{pallit}. Set
$$\norm{f}{Lip(s)}=\supl_k 2^{sk}\norm{P_k f}{L^\infty}.$$
We know that this norm and the norm \eqref{eq.HolNorm} are
equivalent, provided that $0<s<1$. Thus, from these considerations,
we get that $\norm{f}{Lip(s)}$ is bounded whenever $\abs{D}^sf$ is
in the space $L^\infty.$
On the other hand if we choose a weak assumption  $\abs{D}^sf\in L^p,$ for any
$p<\infty$ we get another kind of continuity, more precisely
\beq\label{eq.LebNorm}
 \supl_{x\neq
y}\norm{\frac{\abs{f(x)-f(y)}}{\abs{x-y}^s}}\lp<\infty. 
\eeq
This underlines the fact that the pointwise estimate \eqref{eq.locHol} is more general and could be applied to other cases.
\end{remark}

Lemma \ref{localHo} allows us to prove:

\begin{lem}\label{CommBump}
Let be two function $f, g \in \SSS(\rn)$ and a function $h\in\lo(\rn)$ with $\supp h\subseteq B(0,R).$ Define
\begin{eqnarray}\label{eq.Comm0}
\ch(x)=\int h(y)\round{f(x)-f(x-y)}g(x-y) \, dy,
\end{eqnarray}
then there exists a positive constant $C_{s,n}$, so that the following estimate
\begin{equation}\label{eq.Comm1}
\norm{\ch}\lr=C_{s,n} R^{s}\norm{h}\lo\norm{\abs{D}^s f}\lp \norm{g}\lqq ,
\end{equation}
is satisfied, provided that $1\leq p,q,r \leq \infty, \ \oop+\ooq=\f1r$ and $0<s<1.$
\end{lem}
\bep
We split the proof in two cases, $1< p \leq \infty,$ and $p=1.$\\

\emph{Case $1< p \leq \infty.$ } From the definition of the function $h,$
we can rewrite the function $\ch(x)$ as
\begin{equation}\label{eq.Comm2}
\ch(x)=\intl_{\abs y \leq R} h(y)\round{f(x)-f(x-y)}g(x-y) \, dy,
\end{equation}
and thus, using the the H\"older
inequality,  one has the bound
\begin{eqnarray}\label{eq.Comm3}
 \norm{\ch}\lr &\leq& C \norm{h}\lo \sup_{\abs y \leq R}\round {\int \abs{\round{f(x)-f(x-y)}g(x-y)}^r}^{1/r} \leq
\nonumber\\
 & \leq& C\norm{h}\lo \sup_{\abs y \leq R}\round {\int \abs{f(x)-f(x-y)}^p}^{1/q}\norm{g}\lqq .
\end{eqnarray}

Now, using the pointwise  estimate \eqref{eq.locHol} in the  \eqref{eq.Comm3}, we obtain
\begin{eqnarray}\label{eq.Comm4}
\norm{\ch}\lr\leq C\norm{h}\lo R^{s} \sup_{\abs y \leq R}\norm {M(\abs{D}^s f)(\cdot)+M(\abs{D}^s f)(\cdot-y)}\lp\norm{g}\lqq.
\nonumber\\
\end{eqnarray}
Finally, by the well-known estimate for
the Hardy-Littlewood maximal function $\norm{M(f)}\lp\lesssim \norm{f}\lp,$ with $1<p\leq \infty,$  we obtain  \eqref{eq.Comm1} as desired.\\

\emph{Case $p =1.$ } Consider, for fixed $R\in \rone_{+},$ the function $\Phi\in \lo(\rn)$
given by
\beq\label{eq.L1Func}
\Phi=\Phi_{1}^{R}+\Phi_{2}^{R},
\eeq
where  $\Phi_{1}^{R}$ and $\Phi_{2}^{R}$ are as in \eqref{def.Funct1} and \eqref{def.Funct2}.
From \eqref{eq.RieszPot4},  \eqref{eq.RieszPotFin1} and  \eqref{eq.RieszPotFin2} we have
\begin{eqnarray}\label{eq.Comm1b}
\abs{f(x)-f(x-y)}\leq \abs{I+II}\lesssim \round{\Phi *\abs{\abs{D}^s f(x)}+\Phi*\abs{{\abs{D}^s f(x-y)}}},
\nonumber\\
\end{eqnarray}
with $I$ and $II$ as in the \eqref{eq.RieszPot5a} and \eqref{eq.RieszPot5b}. From the pointwise estimate \eqref{eq.Comm1b} and following now the same line of the proof of the previous case we arrive at
\begin{eqnarray}\label{eq.Comm2b}
\norm{\ch}\lr\leq C\norm{h}\lo  \sup_{\abs y \leq R}\norm {\Phi *\abs{\abs{D}^s f(\cdot)}+ \Phi*\abs{\abs{D}^s f(\cdot-y)}}\lp\norm{g}\lqq.
\nonumber\\
\end{eqnarray}
Finally an application of the Young inequality gives the desired estimate \eqref{eq.Comm1}.

\enp

\subsection{Proof of the Theorem \ref{CommThm}.}
Let us introduce the smooth, positive functions
$$
\chi_j(x) = \vp(2^{-j}x), \ j \geq 1, \  \, \chi_0(x)= \sum_{j\leq 0} \vp(2^{-j}x),
$$
with $\vp$ defined as in \eqref{pallit}, and set
\begin{equation}\label{eq.Projec1}
P=\sum_{j\geq 0}P \chi_j=\sum_{j\geq 0} Q_j,
\end{equation}
where $Q_j=P\chi_j.$ Consider now, for any $k\in\zz,$ the rescaled
functions $\round{Q_j}_{2^{-k}}=2^{nk}Q_j(2^{k}x)$ so that $\supp
\round{Q_j}_{2^{-k}}\subseteq \set{\abs{x}\leq 2^{-k}2^{j}},$ and
they also satisfy:
\begin{equation}\label{eq.Project2}
 \intl \abs{\round{Q_j}_{2^{-k}}} \, dx \lesssim 2^{-j(\ve+s)}.
\end{equation}
To prove \eqref{eq.MainAss} start with the fact that
$$\norm{\int P_k(y)\round{f(x)-f(y)}g(x-y) \, dy}\lr $$
is less than or equal to
 \begin{equation}\label{eq.Project3main}
\suml_{j\geq 0} \norm{\int \round{Q_j}_{2^{-k}}\round{f(x)-f(y)}g(x-y) dy}\lr.
\end{equation}
Because of Lemma \ref{CommBump}, each term of the series
\eqref{eq.Project3main} is less than or equal to a constant times
\begin{eqnarray*}2^{(j-k)s}\norm{\round{Q_j}_{2^{-k}}}\lo\norm{\abs{D}^s f}\lp
\norm{g}\lqq &\lesssim 2^{(j-k)s}2^{-j(\epsilon+s)}\norm{\abs{D}^s
f}\lp \norm{g}\lqq \\ & \sim 2^{-ks}2^{-j\epsilon}\norm{\abs{D}^s
f}\lp \norm{g}\lqq,\end{eqnarray*}and therefore
\eqref{eq.Project3main} is no bigger than a constant times
$2^{-ks}\norm{\abs{D}^s f}\lp \norm{g}\lqq$. This complete the proof
of the the Theorem. \hfill $\Box$

\section{Proof of the embedding Theorem \ref{th:4prima}. }\label{Emb}

First we introduce the following lemma, in which we improve a result
proved in \cite{GeSteTa}. Namely, we have:

\begin{lem} \label{WBS}
If $f$ is a Schwartz function satisfying the property $ {\supp}_{\xi} \, \hat{f}(\xi) \subseteq \{ |\xi| \lesssim 2^k\},$ then
\begin{equation}\label{eq.sum-m}
   \sum_{m } 2^{am} \| [P_k, \vp (2^{-m}\cdot)]
   f\|_{L^2}\leq C \||D|^{-a}f\|_{L^2}
\end{equation}
holds for all integers $k$ and all $ a \in  \round{-n/2, 1}$.
\end{lem}

\bep We shall distinguish two cases.
For $m \leq -k$ we use the estimate \eqref{eq.TW}---with $r=2,
p=\infty, q=2,\ 0<s<1$---to obtain
\begin{equation}\label{eq.estkm0}
    \| [P_k, \vp (2^{-m}\cdot)]
   f\|_{L^2} \lesssim 2^{-ks} \norm{|D|^{s}\vp (2^{-m}\cdot)}\lt  \|f \|_{L^\infty}.
\end{equation}

From Young's inequality we get
$$\norm{|D|^{s}\vp (2^{-m}\cdot)}\lt \leq C 2^{m(n/2-s)};$$

and, from Bernstein's inequality,  we derive
$$ \|f \|_{L^\infty} \lesssim 2^{nk/2} \|f \|_{L^2}.$$
Indeed, from these two estimates we obtain
\begin{eqnarray}\label{eq.estkm}
    2^{am} \| [P_k, \vp (2^{-m}\cdot)]
   f\|_{L^2} \lesssim 2^{m(a+n/2-s)} 2^{k(n/2-s)} \|f \|_{L^2},
\end{eqnarray}
Now, given $a\in(-n/2,1)$, we can find an $s\in (0,1)$ so that
$a+n/2-s>0$. With this $s$ chosen, we get the estimates:
\begin{eqnarray}\label{eq.sum-mneg}
 \sum_{m \leq -k} 2^{am} \| [P_k, \vp (2^{-m}\cdot)]
   f\|_{L^2}&\lesssim&   2^{k(n/2-s)} \sum_{m \leq -k}  2^{m(a+n/2-s)} \|f\|_{L^2}\lesssim
   \nonumber\\
&\lesssim&  2^{-ak} \|f\|_{L^2}.
\end{eqnarray}

For $m > -k$ we follow a similar argument. We replace
\eqref{eq.estkm} by the following commutator estimate
\begin{equation}\label{eq.estkmpog}
    \| [P_k, \vp (2^{-m}\cdot)]
   f\|_{L^2} \lesssim 2^{-ks} \norm{|D|^{s}\vp (2^{-m}\cdot)}{L^\infty} \|f
   \|_{\lt},
\end{equation}
obtained from \eqref{eq.TW}---with $r=2, p=2, q=\infty$. In this
way, since Young's inequality gives
$$\norm{|D|^{s}\vp (2^{-m}\cdot)}{L^\infty} \lesssim 2^{-ms},$$
we arrive at
\begin{eqnarray}\label{eq.summgk}
  \sum_{m > -k} 2^{am} \| [P_k, \vp (2^{-m}\cdot)]
   f\|_{L^2}&\lesssim& 2^{-ks}\sum_{m > -k}  2^{m(a-s)}  \|f\|_{L^2}
   \nonumber\\
 &\lesssim& 2^{-ak} \|f\|_{L^2}
\end{eqnarray}
where we this time we choose $s\in (0,1)$ so that $ a-s <0.$
Finally, this estimate and \eqref{eq.sum-mneg} lead to
\eqref{eq.sum-m}.

\enp

\subsection{Proof of Theorem \ref{th:4prima}.}
Taking into account the definition of $Y^{\ga, d}$ and the mapping
$\phi\to D_x^{-\ga+d/2} \phi=g,$ we need to show that, for every
Schwartz function $g$,
\begin{eqnarray}
\label{eq:970}
\sum_k \round{ \sum_{m} 2^{m/2} \|\vp(2^{-m} \cdot)\tilde{P}_k g\|_{ L^2_x }}^2 \lesssim
\round{\sum_{m} 2^{m/2} \|\vp(2^{-m} \cdot)g\|_{L^2_x}}^2,
\nonumber\\
\end{eqnarray}
 where $\tilde{P}_k$ is the operator defined as in \eqref{eq.tilLP}.
We distinguish two cases. For  $m<-k$ we get
\begin{eqnarray}\label{eq.Besovemb1}
 \sum_k \round{ \sum_{m<-k} 2^{m/2} \|\vp(2^{-m} \cdot)\tilde{P}_k g\|_{ L^2_x }}^2& \leq&
\sum_k   2^{-k} \|\tilde{P}_k g\|_{L_x^2}^2\lesssim
\nonumber \\
& \lesssim& \|D_x^{-1/2}g\|_{L^2_x}
\end{eqnarray}
In the remaining case, $m>-k,$ the decomposition
\begin{eqnarray}\label{eq.dec1}
\vp(2^{-m} x) \tilde{P}_k g=\tilde{P}_k(\vp(2^{-m} \cdot) g)- [\tilde{P}_k, \vp(2^{-m} \cdot)]g.
\end{eqnarray}
is valid. For the first term on the right-hand side of
\eqref{eq.dec1}, we get the estimate
\begin{align}\label{eq.Besovemb2}
&& \round{\suml_k \round{\sum_{m>-k} 2^{m/2}
\|\tilde{P}_k (\vp(2^{-m} \cdot) g\|_{ L_x^2}}^2 }^{1/2}\lesssim
\\
&&\lesssim \sum_m 2^{m/2}
\round{\suml_k \|\tilde{P}_k(\vp(2^{-m} \cdot) g)\|_{ L_x^2}^2 }^{1/2}\sim
\suml_m 2^{m/2} \|  \vp(2^{-m} \cdot)g\|_{ L_x^2}.
\nonumber
\end{align}

It remains to estimate the commutator in \eqref{eq.dec1}. In order
to do that we use the following relation
\begin{eqnarray}\label{eq.dec2}
[\tilde{P}_k, \vp(2^{-m} \cdot)]g=\suml_{l\leq k-N} \tilde{P}_k(\vp(2^{-m} \cdot) P_lg)&+&[\tilde{P}_k, \vp(2^{-m} \cdot)]g_{k-N<\cdot<k+N}+\nonumber\\
&+&\suml_{l\geq k+N} \tilde{P}_k(\vp(2^{-m} \cdot) P_l g),
\end{eqnarray}
with $N\geq 3$, and n general we refer to the first term on the right hand side of the above identity as "high-low interactions", while the second and the third are 
the "low-high" interactions" and ''high-high interaction".  The second term arising above is easy to handle. By
the commutator estimate \eqref{eq.sum-m} (with $a=1/2$), we obtain
\begin{eqnarray}\label{eq.sum-mneg1}
   \sum_{m>-k} 2^{m/2} \| [\tilde{P}_k, \vp(2^{-m} \cdot)]
   g_{k-N<\cdot<k+N}\|_{L^2}&\lesssim&  2^{-k/2}\|g_{k-N<\cdot<k+N}\|_{ L_x^2}.\nonumber\\ 
\end{eqnarray}
Squaring this estimate and taking the sum over $k$ we
find
\begin{equation}\label{eq.l2k}
    \sum_{k} \left( \sum_{m>-k} 2^{m/2} \| \vp (2^{-m}\cdot)
  g_{k-N<\cdot<k+N}\|_{L^2}\right)^2 \lesssim \| D_x^{-1/2}g \|^2_{ L_x^2},
\end{equation}
Now we need to estimate the first and the third terms on the right-hand side of
\eqref{eq.dec2}. In order to do that by the Plancherel identity and the properties of
the operator $P_l$, we obtain the inequalities 
\begin{eqnarray*}
 \| \tilde{P}_k(\vp(2^{-m} \cdot) P_l g)\|_{L_x^2}&\lesssim&\| \tilde{P}_k(\vp(2^{-m} \cdot)\bar{P_l}* P_lg)\|_{L_x^2}\lesssim\\
 &\lesssim& \norm{\int K_{k, m,l} (\xi,\eta) \widehat{P_lg}(\eta) \, d\eta}{
  \lt_\xi}.
\end{eqnarray*}
with $\bar{P_l}=P_{l-1}+P_{l}+P_{l+1}$ and kernel $ K_{k,m,l}(\xi,\eta) $, which satisfies
\begin{equation*}\label{E.6}
|K_{k,m,l}(\xi,\eta)|\leq C 2^{nm} \frac{\vp_k(\xi)\vp_l(\eta)}{2^{M\round{\max(k,l)+m}}}.
\end{equation*}
for some $M>0$ large enough. From these two estimates we get (for more detail on the argument used
here we refer to the paper \cite{GeTa}, Lemmas $5.1$ and $8.1$):
\begin{equation}\label{eq.Kern}
 \| \tilde{P}_k(\vp(2^{-m} \cdot) P_lg)\|_{L_x^2}\lesssim 2^{-M\round{\max(k,l)+m}}\|P_lg\|_{
  \lt_x}.
\end{equation}
for somewhat smaller $M>0.$
We start with the term with indices $k,l$ in the set $\set{l\geq k+N}.$ Let us note that $l\geq k+N\geq -m+N$; that is,
$l+m\geq N$. By the triangle inequality and the estimate \eqref{eq.Kern}, we get, for the
corresponding terms in \eqref{eq.dec2}, the following:
\begin{align}\label{eq.Besovemb3}
& &\suml_k\round {\suml_{m>-k} 2^{m/2}
 \|\suml_{l\geq k+N} \tilde{P}_k(\vp(2^{-m} \cdot) P_lg)\|_{ L^2_x}}^2\lesssim\\
 && \lesssim \suml_k\round{\suml_{m>-k} 2^{m/2}
\suml_{l\geq k+N}  2^{-(l+m)M} \|P_lg\|_{L^2_x}}^2.
\nonumber
\end{align}
Now, taking into account that the chains of inequalities 
\begin{align*}
&&\frac{m}{2}-\round{l+m}M\leq \frac{m}{2}-l-m-\round{l+m}(M-1)\leq \\
&&\leq-\frac{m}{2}-l-\round{l+m}(M-1)\leq -\frac{m}{2}-\frac{l}{2}-\frac{l}{2}-\round{l+m}(M-1)\leq\\
&&\leq -\round{\frac{m}{2}+\frac{l}{2}}-\frac{l}{2}-\round{l+m}(M-1)\leq -\frac{l}{2}-\round{l-k}(M-1),
\end{align*}
are valid, we see that \eqref{eq.Besovemb3} can be controlled by the following
\begin{align}\label{eq.Besovemb3b}
& & \suml_k \round{\suml_{l\geq k+N} 2^{-l/2} 2^{-(l-k)(M-1)}
\|P_lg\|_{ L^2_x}}^2= \\
& &= \sum_{l_1, l_2} 2^{-l_1/2} \|P_{l_1}g\|_{L^2_x}
2^{-l_2/2} \|P_{l_2}g\|_{L^2_x}\suml_{k\leq \min(l_1, l_2)-N}
2^{-(M-1)(l_1+l_2-2k)}.
\nonumber
\end{align}
We claim that if $k\leq \min(l_1, l_2)-N,$ then $l_1+l_2-2k\geq |l_1-l_2|.$ By symmetry, it suffices 
assume $l_2\geq l_1,$ say $l_2=l_1+r$ with $r\geq 0$ and $k= l_1-N-d$ with $d\geq 0.$  From this we get
\begin{align*}
l_1+l_2-2k=2l_1+r-2l_1+2N+2d=r+2N+2d\geq r=\abs{l_1-l_2},
\end{align*}
and by this it is easy to bound \eqref{eq.Besovemb3b} by
\begin{align} \label{eq.Besovemb3c}
& & \sum_{l_1, l_2} 2^{-l_1/2} \|P_{l_1}g\|_{L^2}
2^{-l_2/2} \|P_{l_2}g\|_{L^2_x} 2^{-(M-1)|l_1-l_2|}\lesssim\\
&&\lesssim \suml_{l_1, l_2} 2^{-l_1} \|P_{l_1}g\|_{L^2_x}^2
2^{-(M-1)|l_1-l_2|}\lesssim  \|D_x^{-1/2} g\|^2_{ L^2_x}.
\nonumber
\end{align}
So, from the estimates \eqref{eq.Besovemb3}, \eqref{eq.Besovemb3b} and \eqref{eq.Besovemb3c}
one get the desired
\begin{align}\label{eq.BesovembFin}
\suml_k\round {\suml_{m>-k} 2^{m/2}
 \|\suml_{l\geq k+N} \tilde{P}_k(\vp(2^{-m} \cdot) P_lg)\|_{L^2_x}}^2\lesssim
 \|D_x^{-1/2} g\|^2_{ L^2_x}.
\end{align}

It remains to evaluate the first term on the right hand side of the \eqref{eq.dec2}. We follow the same line of the proof of the previous term, but we consider $l\leq k- N.$ It is enough here to notice
again that $k+m>0$ and use the triangle inequality and the estimate \eqref{eq.Kern}. 
First we notice that in this regime we have
\begin{align*}
&&\frac{m}{2}-\round{k+m}M\leq \frac{m}{2}-\frac{k}{2}-\frac{m}{2}-\round{k+m}\round{M-\frac{1}{2}}\leq \\
&&\leq-\frac{k}{2}-\round{k+m}\round{M-\frac{1}{2}}\leq  -\frac{k}{2},
\end{align*}
and, choosing now $l_2=l_1+r,$ $r\geq 0$ as in the previous case but with  
$k= l_1+r+N+d,$ $d\geq 0,$
\begin{align*}
l_1+l_2-2k=2l_1+r-2l_1-2r-2N-2d\leq r=\abs{l_1-l_2}
\end{align*}
This gives
\begin{align}\label{eq.Besovemb4}
& &\suml_k\round {\suml_{m>-k} 2^{m/2}
 \|\suml_{l\leq k-N} \tilde{P}_k(\vp(2^{-m} \cdot) P_lg)\|_{ L^2_x}}^2\lesssim\\
& & \lesssim \suml_k\round{\suml_{m>-k} 2^{m/2}
\suml_{l\leq k-N}  2^{-(k+m)M} \|P_lg\|_{L^2_x}}^2\lesssim
\nonumber \\
& & \lesssim \suml_k \round{\suml_{l\leq k-N} 2^{-k/2} 
\|P_lg\|_{L^2_x}}^2=
\nonumber \\
& &= \sum_{l_1, l_2} 2^{-l_1/2} \|P_{l_1}g\|_{L^2_x}
2^{-l_2/2} \|P_{l_2}g\|_{L^2_t L^2_x}\suml_{k\geq \max(l_1, l_2)+N}
2^{\frac{1}{2}(l_1+l_2-2k)}
\nonumber\\
& &\lesssim \sum_{l_1, l_2} 2^{-l_1/2} \|P_{l_1}g\|_{L^2}
2^{-l_2/2} \|P_{l_2}g\|_{ L^2_x} 2^{-\frac{1}{2}|l_1-l_2|}\lesssim
\nonumber\\
& & \lesssim \suml_{l_1, l_2} 2^{-l_1} \|P_{l_1}g\|_{L^2_x}^2
2^{-\frac{1}{2}|l_1-l_2|}\lesssim  \|D_x^{-1/2} g\|^2_{L^2_x}.
\nonumber
\end{align}

However, by Sobolev embedding and H\"older's inequality we obtain
\begin{align}\label{eq:971}
&&\|D_x^{-1/2} g\|_{L^2_x}\lesssim \|g\|_{L_x^{\frac{2n}{n+1}}}=\round{\suml_m \intl_{2^{m-1}\leq \abs{x}\leq 2^{m-1}}\abs{g}^{\frac{2n}{n+1}}dx}^{\frac{n+1}{2n}}\lesssim\\
&&\lesssim \round{\suml_m \round{2^{mn}}^{\frac{1}{n+1}} \intl_{2^{m-1}\leq \abs{x}\leq 2^{m-1}}\abs{g}^{\frac{2n}{n+1}}dx}^{\frac{n+1}{2n}}\lesssim
\nonumber\\
&&\lesssim \round{\suml_m 2^{\frac{mn}{n+1}}\norm{\vp(2^{-m} \cdot) g}{L^2_x}}^{\frac{n+1}{2n}}\lesssim
\sum_m 2^{m/2} \|\vp(2^{-m} \cdot) g\|_{\lt_x}\nonumber,
\end{align}
where the last estimate comes from the fact that $n\geq2$ and the property $L^q L^p\subset L^p L^q$ if
$q\leq p.$

Thus, the estimates estimate \eqref{eq.Besovemb1},  \eqref{eq.l2k}, \eqref{eq.BesovembFin}, \eqref{eq.Besovemb3}, \eqref{eq.Besovemb4}  with \eqref{eq:971} and \eqref{eq.Besovemb2} give the proof of
the Theorem \ref{th:4prima}.

\hfill $\Box$

\begin{remark} Note that from the  argument in the proof of the theorem above
one can easily get the following estimate:
\begin{equation}\label{eq.bound}
\norm{P(D)f_k}{Y^{\ga,d}} \lesssim \norm{f_k}{Y^{\ga,d}} =2^{\ga k}
\norm{f}{Y_{d,k}} , \ \forall \, k \in \zz,
\end{equation}
valid for any pseudodifferential operator with symbol $P(\xi) \in
C_0^\infty(\rn).$
\end{remark}

\begin{remark} Some generalizations of the weighted Sobolev norms used in the previous Theorem can
be seen in Theorem 1.6 and Theorem 1.7 in \cite{GeTa}.
\end{remark}

\section{Appendix}
In this section we shall show that the method used to prove the Lemma \ref{CommBump}, case 
$p=1,$ can be applied in order to achieve the general result for $p\in (1, \infty].$ However the use of
the Lemma \ref{localHo} (we referred at this as ''local H\"older continuity'') is of interest in its own right, and such kind has not appeared in the literature. 
So, here we have the following 
\begin{lem}[Version II of Lemma \ref{CommBump}]\label{CommBump2}
Let be $f, g$ two functions in  $\in \SSS(\rn)$ and $h$ a function in $L^1({\rn})$ with $\supp h\subseteq B(0,R).$ Define
\begin{eqnarray}\label{eq.Comm02}
\ch(x)=\int h(y)\round{f(x)-f(x-y)}g(x-y) \, dy,
\end{eqnarray}
then there exists a positive constant $C_{s,n}$, so that the following estimate
\begin{equation}\label{eq.Comm12}
\norm{\ch}\lr=C_{s,n} R^{s}\norm{h}\lo\norm{\abs{D}^s f}\lp \norm{g}\lqq ,
\end{equation}
is satisfied, provided that $1\leq p,q,r \leq \infty, \ \oop+\ooq=\f1r$ and $0<s<1.$
\end{lem}
\bep
We consider now the general case $1\leq p\leq \infty.$\\

Following the line of the proof of the case $p =1,$ we use again the function 
 $\Phi\in \lo(\rn),$ with $R\in \rone_{+}$ fixed, 
given by
\begin{equation*}\label{eq.L1Func2}
\Phi=\Phi_{1}^{R}+\Phi_{2}^{R},
\end{equation*}
where  $\Phi_{1}^{R}$ and $\Phi_{2}^{R}$ are defined in \eqref{def.Funct1} and \eqref{def.Funct2}.
We repeat the proof of the Lemma \ref{localHo}, using \eqref{eq.RieszPot4}, \eqref{eq.RieszPotFin1}  \eqref{eq.RieszPotFin2} (notice that we didn't use the Proposition \ref{Stein} involving the Hardy-
Littlewood maximal operator)
and we arrive at the identity \eqref{eq.Comm1b}. We get easily, by this, the estimate  
\begin{eqnarray}\label{eq.Comm2b2}
\norm{\ch}\lr\leq C\norm{h}\lo  \sup_{\abs y \leq R}\norm {\Phi *\abs{\abs{D}^s f(\cdot)}+ \Phi*\abs{\abs{D}^s f(\cdot-y)}}\lp\norm{g}\lqq.
\nonumber\\
\end{eqnarray}
Now, by the Young inequality again, we bound the sum on the left hand side of the above estimate \eqref{eq.Comm2b2}
in the following way
\begin{eqnarray}\label{eq.Comm2bBound}
 \sup_{\abs y \leq R}\norm {\Phi *\abs{\abs{D}^s f(\cdot)}+ \Phi*\abs{\abs{D}^s f(\cdot-y)}}\lp\lesssim
 R^s \norm{\abs{D}^s f}\lp.
\nonumber\\
\end{eqnarray}
This estimate and the \eqref{eq.Comm2b2} give the estimate \eqref{eq.Comm12}. The alternative proof is now complete.

\enp

\end{document}